\documentclass{amsart}

\usepackage{amsmath}
\usepackage{amssymb}
\usepackage{amsfonts}
\usepackage{graphicx}
\usepackage{texdraw}
\usepackage{tikz}
\usepackage{enumitem}
\usepackage{graphpap}
\usepackage{wasysym}
\usepackage{color}
\usepackage{pxfonts}
\usepackage[OT2,T1]{fontenc}

\theoremstyle{plain}

\theoremstyle{remark}

\newtheorem*{rem}{Remark}

\newcommand{\C}{{\mathbb C}}
\newcommand{\E}{{\mathbb E}}

\newcommand{\eps}{{\varepsilon}}

\renewcommand{\d}{{\partial}}
\newcommand{\dbar}{\bar{\partial}}
\newcommand{\1}{\mathbf{1}}

\begin{document}

\title[]{A word on thermal equilibrium and fluctuations of two-dimensional Coulomb systems}


\author{Yacin Ameur}

\begin{abstract}  In this note we examine the relationship between two-dimensional Coulomb systems and the thermal equilibrium measure, such as defined in the lecture notes \cite{Se}, pointing out some of its discrepancies with respect to the one-particle distribution. 
\end{abstract}

\maketitle


For the benefit of mutual understanding we will here synchronize notation with respect to a few notions such as ``expectation of fluctuations'' and ``thermal equilibrium measure''.

It is perhaps surprising, but when this is done, 
it  looks like some inconsistencies emerge, which might require some further clarification. 

Regardless, we hope that the following remarks will ultimately contribute to dispersing some apparent conceptual confusion in the area.


\smallskip


So let us start. 
In \cite[eq (1.1.1)]{Se} (with $s=0$, $d=2$) the Hamiltonian is defined as
$$H_n=\frac 1 2 \sum_{j\ne k}\log \frac 1 {|z_j-z_k|}+n\sum_{j=1}^n V(z_j).$$

The density of the Gibbs measure (1.1.5) in \cite{Se} is proportional to
\begin{equation}\label{exp}\exp(-\beta  H_n(z_1,\ldots,z_n)).\end{equation}

To compare with \cite{AT} (as well as with \cite{BF}, cf.~also e.g.~\cite{B}) we now write
$$V=Q/2.$$
So the exponent in \eqref{exp} is the negative of
$$\frac \beta {2}(\sum_{j\ne k}\log\frac 1 {|z_j-z_k|}+n\sum_{j=1}^n Q(z_j)),$$
which is our form (in \cite{AT}) of writing the exponent of the Gibbs measure but with $\beta$ replaced by $\beta/2.$

Next consider the energy-entropy functional from \cite[eq (2.0.3)]{Se}. Let us call it $\tilde{F}(\mu)$. It is
\begin{equation}\label{sdef}\tilde{F}(\mu)=\tilde{I}(\mu)+\frac 1 \theta E(\mu)\end{equation}
where, by \cite[eq (2.0.2)]{Se},
$$\tilde{I}(\mu)=\frac 1 2 \int_{\C^2}\log \frac 1 {|z-w|}\, d\mu(z)d\mu(w)+\int_\C V\, d\mu,$$
so, according to our notation form \cite{AT},
$$\tilde{I}(\mu)=\frac 1 2 I_Q(\mu);\qquad (I_Q(\mu):=\int_{\C^2}\log \frac 1 {|z-w|}\, d\mu(z)d\mu(w)+\int_\C Q\, d\mu).$$

In both papers, the entropy is defined as
$$E(\mu)=\int_\C \mu\log \mu.$$

Finally, according to Remark 2.21 in \cite{Se}, we should take
$$\theta=n\beta,$$
so we have
$$\tilde{F}(\mu)=\frac 1 2 I_Q(\mu)+\frac 1 {n\beta}E(\mu)=\frac 1 2 (I_Q(\mu)+\frac 2 {n\beta}E(\mu)).$$

If we put $\beta_*=\beta/2$ and $F(\mu)=I_Q(\mu)+\frac 1 {n\beta_*}E(\mu)$ (as in \cite{AT} but with $\beta_*$) we see that
$$\tilde{F}(\mu)=\frac 1 2 F(\mu).$$

The last relationship shows that the minimizers of $\tilde{F}$ and of $F$ are the same, i.e., we obtain the thermal equilibrium measure in the determinantal case by minimizing $F(\mu)$ where $\beta_*=1$.


Furthermore we note that the formula for the bulk density in the determinantal case written in \cite[Remark 2.21]{Se} is taken from \cite{BF} and is incorrectly reproduced. 

According to our conventions (see \cite[eq (6.14)]{AT} or \cite[eq (5.24)]{BF}) the bulk density is proportional to
$$\Delta Q+n^{-1}\frac 1 2 \Delta\log\Delta Q+O(n^{-2})$$
where $\Delta$ is, say, the standard Laplacian.

Setting $Q=2V$ this becomes
$$2\Delta V+n^{-1}\frac 1 2 \Delta\log\Delta V+O(n^{-2}),$$
or, dividing by $2$,
\begin{equation}\label{c1}\Delta V+n^{-1}\frac 1 4 \Delta\log\Delta  V+O(n^{-2}).\end{equation}

However \cite[eq (2.5.32)]{Se} states that the thermal equilibrium measure (for $\beta=2$) is proportional to
\begin{equation}\label{c2}\Delta V+n^{-1}\frac 1 2 \Delta\log\Delta V+O(n^{-2}).\end{equation}

As in \cite[Theorem 6.2]{AT}, we confirm that \eqref{c1} and \eqref{c2} are inconsistent with the hypothesis
the coefficient for $n^{-1}$ is the same for the
thermal equilibrium measure and for the one-point function. 

As remarked in \cite{AT}, the equation \eqref{c2} is also inconsistent with the limiting expectation of fluctuations, a distribution called $f\mapsto \rho_{1/2}(f)$ in \cite{A,ZW} (see eq (5.16)), also known as
$f\mapsto e_f$ in eg \cite{AC0}. 

Indeed, writing $\beta_*=\beta/2$ we have by \cite[(5.16)]{ZW} that (with $L(z)=\log\Delta Q(z)$ in a neighbourhood of $S$ and $L^S$ the usual Poisson modification in the exterior)
\begin{align*}\rho_{1/2}(f)&=(\frac 2 {\beta_*}-1)\int_S f\cdot \frac 1 2 \Delta L\, dA+(\frac 2 {\beta_*}-1)\frac 1 {8\pi}\int_{\d S}f\cdot \d_{\tt{n}}(L^S-L)\, ds\\
&+(\frac 2 {\beta_*}-1)\frac 1 {8\pi}\int_{\d S}\d_{\tt{n}}f\, ds.
\end{align*}

(Here and from now on we write $\Delta:=\d\dbar$ and $dA(z):=\frac 1 \pi d^2z$; the normal derivative is taken in the exterior of the droplet.)

The leading term of the suitably normalized density $\rho$ is the equilibrium measure $\rho_0=\Delta Q\cdot \1_S$, and we have a large $n$ expansion (with distributional coefficients)
$$\rho=\rho_0+n^{-1}\rho_{1/2}+n^{-2}\rho_1+\cdots.$$

In particular, the bulk density is (in a distributional sense) proportional to
\begin{equation}\label{heur}\Delta Q+n^{-1}(\frac 1 {\beta_*}-\frac 12)\Delta\log\Delta Q+O(n^{-2}),\end{equation} whereas the thermal equilibrium density in the bulk has the asymptotic
$$\Delta Q+n^{-1}\frac 1 {\beta_*}\Delta\log \Delta Q+O(n^{-2}).$$ 
We now show that (at least for $\beta_*=1$) the asymptotics of the thermal equilibrium also differs in an essential way from $\rho_{1/2}$ near the edge.

To this end we consider the minimizer $\mu_{1,n}$ of the modified functional
$$F_1(\mu)=I_Q(\mu)+\frac 1 {n\theta_1}E(\mu),\qquad \frac 1 \theta_1=\frac 1 {\beta_*}-\frac 1 2.$$
It might be hoped that $\mu_{1,n}$ is a good approximation to the density, since in the bulk the subleading terms coincide in the sense of distributions. However, when $\beta_*=1$, the density of $\mu_{1,n}$ differs from the $1$-point function in the leading order at the boundary. 
This is easily seen in the case of the Ginibre ensemble as follows. 

As in \cite[Sec 6]{AT} we see that the density $\delta_{1,n}$, ($d\mu_{1,n}=\delta_{1,n}\, dA$), obeys the variational equation
\begin{equation}\label{var}-\delta_{1,n}+\Delta Q+\frac 1 {n\theta_1}\Delta \log \delta_{1,n}=0.\end{equation}

We take $Q(z)=|z|^2$ so that $S=\{|z|\le 1\}$ and we rescale about the boundary point
$z=1$,
$$z=1+\frac 1 {\sqrt{n}}u,\qquad \tilde{\delta}_{1,n}(u):=\delta_{1,n}(z).$$

Similarly we rescale the $1$-point intensity function (density with respect to $dA$)
$$R_n(z)=\lim_{\eps\to 0}\frac {\E_n(\#\{z_j\}_1^n\cap D(z,\eps))}{\eps^2}$$
(which in our normalization roughly approximates $n\Delta Q\1_S$) and form
$$g_{n}(u)=\frac 1 n R_n(z).$$

A computation in \cite[p. 24]{AT} yields that 
$$\lim_{n\to\infty}\Delta \log g_n(0)=-\frac 2 \pi.$$
Also $g_n(0)\to 1/2$ as $n\to\infty$.

If we assume that $\tilde{\delta}_{1,n}(0)\to 1/2$ (along some subsequence) we now deduce a contradiction. Indeed by rescaling in \eqref{var} we find
$$-\tilde{\delta}_{1,n}(u)+1+\frac 1 {\theta_1}\Delta\log \tilde{\delta}_{1,n}(u)=0.$$
Setting $\beta_*=1$ we have $\theta_1=2$. Next setting $u=0$ and using $\tilde{\delta}_{1,n}(0)\to 1/2$ we conclude that $\Delta\log \tilde{\delta}_{1,n}(0)\to -1$ as $n\to\infty$. Since $2/\pi\ne 1$, we see that the approximation of $R_n$ by $n\delta_{1,n}$ fails at the boundary point, and actually $\liminf_{n\to\infty}\{n^{-1}\|R_n-n\delta_{n,1}\|_\infty\}>0$.

By similar computations we easily see that there is no way to tweak the parameter $\theta$ in \eqref{sdef} so that the resulting thermal equilibrium measure coincides with the distribution $\rho(f)$ beyond the leading term (equilibrium density) $\rho_0$. 

We remark that the correction term $\rho_{1/2}$ (also for perturbations of the potential) is central in \cite{AC0,AHM} and related works, where it is used to deduce Gaussian field convergence, among other things. (Cf.~\cite{A} for another recent application.


\begin{rem}
The factor $(\frac 1 {\beta_*}-\frac 1 2 )$ appearing in the functional $\rho_{1/2}$ is very natural, eg from the point of view of Ward identities. The term $\frac 1 {\beta_*}$ comes from changing variables in the energy expression whereas the $-\frac 1 2$ comes from the corresponding change of integration measure. See also \cite{ZW,Z}, for instance.
\end{rem}



\begin{thebibliography}{999}
\bibitem{A} Ameur, Y., \textit{A formula for the edge density $\sqrt{n}$-correction for two-dimensional Coulomb systems,} , arxiv: 2510.16945.
\bibitem{AC0} Ameur, Y., Cronvall, J., \textit{On fluctuations of Coulomb systems and universality of the Heine distribution}, arxiv: 2411.10288.
\bibitem{AHM} Ameur, Y., Hedenmalm, H., Makarov, N., \textit{Random normal matrices and Ward identities}, Ann. Probab. \textbf{43} (2015), 1157--1201.
\bibitem{AT} Ameur, Y., Troedsson, E., \textit{Remarks on the one-point density of Hele-Shaw $\beta$-ensembles}, Rev. Mat. Iberoam. \textbf{41} (2025), 1555--1582.
\bibitem{AS} Armstrong, S., Serfaty, S., \textit{Thermal approximation of the equilibrium measure and obstacle problem}, Ann. Fac. Sci. Toulouse Math. (6) \textbf{31} (2022), no. 4, 1085–1110.
\bibitem{B} Berman, R., \textit{Bergman kernels and weighted equilibrium measures in $\C^n$}. Indiana Univ. Math. J. \textbf{58}
(2009), 1921-1946.
\bibitem{BF} Byun, S.-S., Forrester, P.J., \textit{Progress on the study of the Ginibre ensembles I: GinUE,} KIAS Springer Series in Mathematics \textbf{3} (2025).
\bibitem{Se} Serfaty, S., \textit{Lecture Notes on Coulomb and Riesz gases}, arXiv:2407.21194v1.
\bibitem{ZW} Zabrodin, A., Wiegmann, P., \textit{Large $N$ expansion for the 2D Dyson gas}, J. Phys. A \textbf{39} (2006), no.28, 8933--8964.
\bibitem{Z} Zabrodin, A., \textit{Matrix models and growth processes: from viscous flows to the quantum Hall effect}, arxiv: 0412219.
\end{thebibliography}
\end{document}